\documentclass[11pt]{amsart}
\usepackage{amssymb,amscd}

\numberwithin{equation}{section}

\newtheorem {Theorem}		{Theorem}
\newtheorem {T}	{Th\'eor\`eme}
\newtheorem*{W}{Weinstein's Theorem}
\newtheorem*{w}{Th\'eor\`eme de Weinstein}
\newtheorem* {Corollary} {Corollary}
\newtheorem* {Corollaire} {Corollaire}

\theoremstyle{remark}

\def \rpo	{relative periodic orbit}
\def \re	{relative equilibrium}
\def \augh      {h - \langle\Phi | \xi\rangle}
\def \th        {{\tilde h}}
\def \inv	{^{-1}}

\newcommand{\fg} {{\mathfrak g}}
\newcommand{\fh} {{\mathfrak h}}
\newcommand{\fl} {{\mathfrak l}}

\font\boo=cmr5
\font\baa=cmbxti10

\begin{document}
\title{On relative normal modes}

\author{E.~Lerman}\author{T.~F.~Tokieda}

\address{Department of Mathematics, University of Illinois, 
Urbana IL 61801, USA}
\email{lerman@math.uiuc.edu}\email{tokieda@math.uiuc.edu} 

\date{28/10/98}

\begin{abstract}
	We generalize the Weinstein-Moser theorem on the existence of 
nonlinear normal modes near an equilibrium in a Hamiltonian system to 
a theorem on the existence of \rpo s near a \re\
in a Hamiltonian system with continuous symmetries.  In
particular we prove that under appropriate hypotheses there exist
\rpo s near relative equilibria even when these
relative equilibria are singular points of the corresponding moment
map, i.e. when the reduced spaces are singular.

\bigskip

\bigskip

\noindent {\bf Sur les oscillations normales relatives}

\bigskip

\noindent {R{\boo \'ESUM\'E}.}~~ On g\'en\'eralise le th\'eor\`eme 
d'existence de
Weinstein-Moser des oscillations normales non-lin\'eaires voisines
d'un \'equilibre dans un syst\`eme hamiltonien, \`a un th\'eor\`eme
d'existence des orbites p\'eriodiques relatives voisines d'un
\'equilibre relatif dans un syst\`eme hamiltonien \`a sym\'etrie 
continue.  Entre autres on d\'emontre, moyennant quelques
hypoth\`eses, qu'il existe des orbites p\'eriodiques relatives
pr\`es des \'equilibres relatifs m\^eme lorsque ceux-ci sont des
points singuliers de l'application de moment, c'est-\`a-dire lorsque
les espaces r\'eduits sont singuliers.
\end{abstract}

\maketitle

\noindent {\baa Version fran{\c c}aise abr\'eg\'ee}.

\smallskip

	Dans cette note on g\'en\'eralise le th\'eor\`eme de
Weinstein-Moser \cite{W1, Ms, W2, MnRS, Ba} sur les oscillations
normales non-lin\'eaires voisines d'un \'equilibre dans un syst\`eme
hamiltonien, \`a un th\'eor\`eme sur les orbites p\'eriodiques
relatives voisines d'un \'equilibre relatif dans un syst\`eme 
hamiltonien sym\'etrique.

	Soit $M$ une vari\'et\'e symplectique munie d'une op\'eration
hamiltonienne d'un groupe de Lie connexe $G$ et donc d'une application
de moment $\Phi: M \to \fg^{\ast}$.  Une orbite du champs hamiltonien
$X_h$ d'un hamiltonien $G$-invariant $h\in C^{\infty}(M)^G$ est un
{\it \'equilibre relatif\/} [resp.\ une {\it orbite p\'eriodique
relative\/}] si son image dans $M/G$ est un point [resp.\ un lacet].

	Les am\'eliorations ult\'erieures du th\'eor\`eme
de Weinstein-Moser se pr\^etent toutes \`a des g\'en\'eralisations
par notre m\'ethode; n\'eanmoins, afin de ne pas alourdir l'expos\'e
on traitera seulement la version originale:

\begin{w}\label{Theorem1} {\rm [W1]}
	Soit $h$ un hamiltonien sur un vectoriel symplectique $V$
avec $dh(0)$ nul et $d^2h(0)$ d\'efini positif.  Alors 
pour tout $\varepsilon > 0$ petit, $h\inv(h(0) + \varepsilon)$
porte au moins $\frac{1}{2}\dim V$ orbites p\'eriodiques de $X_h$.
\end{w}

	Soit $x\in M$ un point d'un \'equilibre relatif pour
un hamiltonien $G$-invariant.  Existe-t-il des orbites p\'eriodiques
relatives voisines?  A cela nulle difficult\'e si $x$ est un point 
r\'egulier de $\Phi$, ou m\^eme dans le cas singulier si la dimension
de la strate de son image $\bar{x}$ dans l'espace r\'eduit est
strictement positive, car le probl\`eme se ram\`ene alors au th\'eor\`eme
de Weinstein sur cette strate, qui est stable par rapport \`a la
dynamique de l'espace r\'eduit stratifi\'e \cite{AMM, SL}.  Mais
que dire si la strate de $\bar{x}$ est un point?

	On verra qu'en l'occurrence il existe encore des
familles d'orbites p\'eriodiques relatives voisines de l'\'equilibre
relatif, pourvu qu'un certain `rempla{\c c}ant' du hessien soit
d\'efini en tant que forme quadratique 
et que le groupe d'isotropie de $\Phi(x)$ soit un tore. 
On a en effet un th\'eor\`eme de correspondance que voici:

\begin{T}\label{Theorem3}
	Soit $M$ une vari\'et\'e symplectique munie d'une op\'eration
hamiltonienne d'un groupe de Lie connexe $G$ et d'une application
de moment $\Phi: M \to \fg^{\ast}$, et soit $h\in C^{\infty}(M)^G$ un
hamiltonien $G$-invariant. 

	Supposons que $x\in M$ appartienne \`a un \'equilibre relatif
pour $h$ et que le groupe d'isotropie de $\Phi(x)$ soit un tore.
Alors il existe un vectoriel symplectique $V$ muni d'une
op\'eration hamiltonienne du groupe d'isotropie $G_x$ de $x$ et
un hamiltonien $G_x$-invariant $h_V\in C^{\infty}(V)^{G_x}$ tels que
\begin{enumerate}
	\item l'origine $0\in V$ est un \'equilibre pour $h_V$;

	\item le hessien $d^2h_V(0)$ de $h_V$ peut \^etre calcul\'e
	\`a partir de $h$;

	\item toute orbite p\'eriodique $G_x$-relative pour $h_V$
	dans $V$ donne lieu \`a une famille d'orbites p\'eriodiques
	$G$-relatives pour $h$ dans $M$.
\end{enumerate}
\end{T}

	Ici $V$ d\'esigne le sous-espace symplectique maximal de
$\ker d\Phi(x)$ (dit {\it tranche symplectique\/} en $x$) et 
$d^2h_V(0)$ est la forme quadratique induite sur $V$ par 
$d^2(\augh)(x)|_{\ker d\Phi(x)}$, o\`u
$\xi$ est une {\it vitesse\/} de $x$ que l'on trouve dans
$$
x \text{ appartient \`a un \'equilibre relatif pour } h \ 
\Longleftrightarrow\ \exists \xi\in \fg, \quad d(\augh)(x) = 0.
$$
La forme normale locale de Guillemin-Sternberg et de Marle
\cite{GS, Mr} permet de construire une fonction $h_V$ sur $V$
dont le hessien est $d^2h_V(0)$.  Cette notion de hessien 
intervient dans l'\'etude de la 
stabilit\'e et des bifurcations des \'equilibres relatifs aux 
points singuliers de l'application de moment \cite{LS, Ortega, OR}.

  	Le Th\'eor\`eme 1, joint par exemple au th\'eor\`eme de Weinstein,
conduit au r\'esultat escompt\'e:

\begin{Corollaire}\label{Corollary}
	Soient $M$, $G$, $\Phi : M \to \fg^{\ast}$, 
$h\in C^{\infty}(M)^G$ comme dans le Th\'eor\`eme 1.  Supposons
que $x\in M$ appartienne \`a un \'equilibre relatif et que le groupe
d'isotropie de $\Phi(x)$ soit un tore; appelons $V$ la tranche symplectique
en $x$.  S'il existe une vitesse $\xi$ de $x$ telle que 
$d^2(\augh)(x)|_V$ est d\'efini positif, alors pour tout 
$\varepsilon > 0$ petit, 
$\{ y \in M \mid (\augh)(y) = (\augh)(x) + \varepsilon \}$
porte des familles d'orbites p\'eriodiques relatives.
\end{Corollaire}

    	On peut sans doute se passer de l'hypoth\`ese selon laquelle
le groupe d'isotropie de $\Phi(x)$ est un tore, mais la d\'emonstration
ne promet gu\`ere d'\^etre si simple.  Signalons que dans
le cas o\`u $G$ est compact,
l'hypoth\`ese est satisfaite de fa{\c c}on {\it g\'en\'erique}.

\smallskip

\noindent {\it Esquisse de la d\'emonstration.}

\smallskip  

	D'abord, on se ram\`eme au cas o\`u $G$ est un tore en passant
\`a un sous-syst\`eme hamiltonien dont le groupe de sym\'etrie est le
groupe d'isotropie de $\Phi(x)$, gr\^ace \`a un th\'eor\`eme de
Guillemin-Sternberg (voir \cite{GLS}, Corollaire 2.3.6).  Ensuite,
on plonge, par un symplectomorphisme \'equivariant, un voisinage
de $G\cdot x$ dans $(T^{\ast}L\times V)/\Gamma$ ($L$ est le tore
compl\'ementaire de la composante connexe $K$ de $1$ dans $G_x$, et 
$\Gamma = G_x/K$); en r\'eduisant par $L$, on r\'eduit ce dernier
espace \`a $V$ et $h$ \`a $h_V$.  Enfin, on v\'erifie que les orbites
p\'eriodiques relatives pour $h_V$ se rel\`event \`a des orbites
p\'eriodiques relatives pour $h$: les unes et les autres correspondent
aux orbites p\'eriodiques du syst\`eme obtenu soit en r\'eduisant
$V$ par $\Gamma\times K$, soit en r\'eduisant 
$(T^{\ast}L\times V)/\Gamma$ par $L\times K$.
 
\smallskip

	Comme un exemple d'application, on remarque que le Corollaire 
assure l'existence des orbites p\'eriodiques relatives d'une paire
de corps rigides \`a sym\'etrie axiale li\'es par un potentiel
qui d\'epend de l'angle entre les corps.

\bigskip \bigskip

\hrule	
	
\bigskip \bigskip

\section{Introduction.}

In this paper we generalize the Weinstein-Moser theorem (\cite{W1, Ms,
W2, MnRS, Ba} and references therein) on the nonlinear normal modes
near an equilibrium in a Hamiltonian system to a theorem on the
\rpo s near a \re\ in a symmetric Hamiltonian system.

Let $M$ be a symplectic manifold, with a Hamiltonian action
of a connected Lie group $G$ and a moment map $\Phi: M \to
\fg^{\ast}$.  Recall that an orbit of the Hamiltonian vector field
$X_h$ of a $G$-invariant Hamiltonian $h\in C^{\infty}(M)^G$ is a {\it
\re\/} if its image in the orbit space $M/G$ is a
point, and that an orbit of $X_h$ is a {\it \rpo\/} if its image in
$M/G$ is a closed curve.

	Later improvements of the Weinstein-Moser theorem lend
themselves to generalizations by our method; to streamline the
presentation, however, we shall treat only the original version:

\begin{W}\label{Theorem1} {\rm [W1]}
Let $h$ be a Hamiltonian on a symplectic vector space $V$ such that
its differential at the origin $dh(0)$ is zero and its Hessian at the
origin $d^2h(0)$ is positive definite.  Then for every small
$\varepsilon > 0$, the energy level $h\inv(h(0) + \varepsilon)$
carries at least $\frac{1}{2}\dim V$ periodic orbits of the Hamiltonian
vector field of $h$.
\end{W}

Now suppose $x\in M$ lies on a \re\ for a
$G$-invariant Hamiltonian $h$. If $x$ is a {\it regular\/} point of
the moment map $\Phi$, then the reduced space at $\mu = \Phi(x)$ is 
smooth near the image $\bar{x}$ of $x$. 
Provided appropriate conditions hold on the Hessian of the reduced
Hamiltonian, Weinstein's theorem applied to the reduced system
guarantees at least $\frac{1}{2} \dim (\text{reduced
space})$ periodic orbits, and hence as many families of \rpo s 
near the \re.  (If $x$ lies on a
\rpo, then the orbit through $g\cdot x$ is relative periodic also
for every $g\in G$.)\ \ On the other hand, if $x$ is
a {\it singular\/} point of $\Phi$, then the reduced space at
$\mu$ is a stratified space, and the reduced dynamics preserves
the stratification \cite{AMM, SL}.  Unless the stratum
through $\bar{x}$ is an isolated point, we have again
$\frac{1}{2}\dim (\text{stratum})$ families of \rpo s, provided 
appropriate conditions hold on the Hessian of the
restriction of the reduced Hamiltonian to the stratum.  But what if
the stratum through $\bar{x}$ is an isolated point? It is difficult to
make sense of Hessians on singular spaces.

	 We shall show that in this
case also there are families of \rpo s near the \re\
provided a certain substitute for the Hessian is
definite and the isotropy group $G_{\mu}$ of $\mu$ is a torus.
The proof amounts to a computation in `good coordinates'  that allows 
us to reduce our problem to Weinstein's theorem.  We proceed via
the following `correspondence theorem':

\begin{Theorem}\label{Theorem2}
       	Let $M$ be a symplectic manifold, with a Hamiltonian
action of a connected Lie group $G$ and a moment map $\Phi: M \to
\fg^{\ast}$, and let $h\in C^{\infty}(M)^G$ be a $G$-invariant Hamiltonian.

Suppose $x\in M$ lies on a \re\ for $h$ and the isotropy group of 
$\Phi(x)$ is a torus.  Then there exist a
symplectic vector space $V$ with a Hamiltonian action of the
isotropy group $G_x$ of $x$ and a $G_x$-invariant Hamiltonian $h_V\in
C^{\infty}(V)^{G_x}$ such that
\begin{enumerate}
	\item the origin $0\in V$ is an equilibrium for $h_V$;

	\item the Hessian $d^2h_V(0)$ of $h_V$ can be computed in
	terms of $h$;

	\item every $G_x$-\rpo\ for $h_V$ on $V$ sufficiently close 
	      to the origin gives rise to a family of $G$-\rpo s
	      for $h$ on $M$.
\end{enumerate}
\end{Theorem}

	The meaning of the vector space $V$ and of the quadratic form
$d^2h_V(0)$ is as follows.  Note that
$$
x\in M \text{ lies on a \re\ for } h \ \Longleftrightarrow\
d(\augh)(x) = 0 \; \text{ for some } \xi\in \fg.
$$
The vector $\xi$ is called a {\em velocity\/} of $x$.  Velocity is 
not unique: any two velocities of $x$ differ by a vector in the 
isotropy Lie algebra $\fg _x$ of $x$.

	The function $\augh$ is constant along 
the orbit $G_{\mu}\cdot x$, where as above $G_{\mu}$ is the isotropy 
group of $\mu = \Phi (x)$.  The form $d^2(\augh)(x)|_{\ker d\Phi(x)}$ 
is therefore always degenerate and descends
to a well-defined form on the vector space $V := \ker d\Phi (x) / T_x
(G_{\mu}\cdot x)$; alternatively we can think of $V$ as the maximal
symplectic subspace of $\ker d\Phi (x)$.  $V$ is called the {\em
symplectic slice\/} at $x$.  It follows from the local normal form
theorem of Guillemin-Sternberg and Marle \cite{GS, Mr} that there exists a 
$G_x$-invariant symplectic submanifold $\Sigma$ passing through $x$
such that
$T_x \Sigma = V$.  Thus, locally $\Sigma \simeq V$ as symplectic manifolds,
with $x$ corresponding to the origin in $V$.  The function $h_V$ in
Theorem~\ref{Theorem2} is the restriction $(\augh)|_\Sigma =(\augh)|_V $.  
Since Hessians are functorial, 
$d^2 h_V (0) = d^2(\augh)|_V(0) = d^2(\augh)(x) |_V.$ 
This notion of Hessian has been used in \cite{LS, Ortega, OR} to study
the stability and bifurcation of relative equilibria at singular points
of the moment map.

As an example of applications of Theorem~\ref{Theorem2}, we combine
it with Weinstein's theorem to obtain:

\begin{Corollary}\label{Corollary}
Let $M$, $G$, $\Phi :M \to \fg^* $, $h\in C^\infty (M)^G$
be as in Theorem~\ref{Theorem2}.  Suppose $x\in M$ lies on a
\re\ for $h$ and the isotropy group of $\Phi(x)$ is a torus;
call $V$ the symplectic slice at $x$. 
If there is a velocity $\xi$ of $x$ for which $d^2(\augh)(x)|_V$ 
is positive definite, then for every small $\varepsilon > 0$, 
the level set $\{ y \in M \mid (\augh)(y) = (\augh)(x) + \varepsilon\}$
carries families of \rpo s.
\end{Corollary}

	We expect that the assumption on the isotropy group of
$\Phi(x)$ being a torus can be dropped, but the proof is unlikely to
be quite so simple.   In case $G$ is compact, the assumption is 
satisfied {\it generically}.

\section{ \protect Proof of Theorem~\ref{Theorem2}.}

	Let $M$, $G$, $\Phi: M \to \fg^{\ast}$, $h\in
C^{\infty}(M)^G$ be as in the statement of Theorem~\ref{Theorem2}. 
We are supposing that $x\in M$ lies on a \re\ for $h$ and that the 
isotropy group $G_{\mu}$ of $\mu = {\Phi(x)}$ is a torus.

	First, we show that it suffices to consider the case when the
whole $G$ is a torus.  Since $G_{\mu}$ is compact, we can produce an 
$Ad^\dagger (G_{\mu})$-invariant inner
product on $\fg^*$ and take the corresponding $G_{\mu}$-equivariant 
splitting $\fg = \fg_{\mu} \oplus \fh$.  Then a small enough 
$G_{\mu}$-invariant
neighborhood $B$ of $\mu$ in the affine plane $\mu + \fh^\circ $ is
transverse to the moment map (here $\fh^\circ $ denotes the
annihilator of $\fh $ in $\fg$).  Hence $R := \Phi \inv (B)$ is a
$G_{\mu}$-invariant submanifold of $M$ containing $x$.  In fact,
by the symplectic cross-section theorem of Guillemin and Sternberg
(cf.\ \cite{GLS}, Corollary~2.3.6), $R$ is a
{\em symplectic\/} submanifold of $M$ and the action of $G_{\mu}$ on $R$ is
Hamiltonian; the moment map for this action is the restriction
of $\Phi$ to $R$ followed by the natural projection $\fg^* \to
\fg_{\mu}^*$.  Since $\fg^* \to \fg_{\mu}^*$ restricted to 
$\mu + \fh^\circ$ is an isomorphism, the moment map for the
action of $G_{\mu}$ on $R$ is $\Phi|_R$ up to a linear isomorphism.  It
follows that 
$$
	\ker d\Phi(y) = \ker d(\Phi|_R)(y)
$$
for every $y\in R$ (cf.\ \cite{LS}, Lemma~2.5).
Moreover, because $h$ is $G$-invariant, the flow of $X_h$
preserves the fibers of the moment map, and so the flow
preserves $R$.  It follows that
$$
	(X_h)|_R = X_{(h|_R)}
$$
(cf.\ \cite{L}, p.218).  Thus, we have found a Hamiltonian sub-system
$(R, G_{\mu}, \Phi|_R, h|_R)$ for which the symmetry group $G_{\mu}$ 
is a torus.  Passing to this sub-system, we may and shall assume 
without loss of generality that $G$ is a torus and $G = G_{\mu}$.

	Second, we construct the Hamiltonian $h_V$.  The connected
component $K$ of $1$ in the isotropy group $G_x$ of $x$ is a torus and
has a complementary torus $L$, so that $G = L\times K$.  The finite
abelian group $\Gamma = G_x/K$ may be regarded as a subgroup of $L$.
Then $\Gamma$ acts symplectically (by
the lifted action) on $T^*L = L\times \fl^*$ and on the symplectic
slice $V$ at $x$.  Hence $\Gamma$ acts diagonally on $T^*L \times V$.  Note
that $G$ too acts on $T^*L \times V$: $L$ acts on its own cotangent 
bundle and $K$ acts on $V$.  Hence $G$ acts in a Hamiltonian
fashion on $(T^*L \times V)/\Gamma$.  The orbit of $[1,0,0] \in (T^*L
\times V)/\Gamma$ is isotropic and is isomorphic to $G/(\Gamma\times K)
\simeq G\cdot x$.  By the equivariant isotropic
embedding theorem, there exist $G$-invariant neighborhoods $U_0$ of
$[1,0,0]$ in $(T^*L \times V)/\Gamma$ and $U$ of $x$ in $M$ and an
equivariant symplectomorphism $\sigma : U_0\to U$ such that $\sigma
([1,0,0]) = x$ and the derivative $d \sigma ([1,0,0]) $ sends
$V\subset T_{[1,0,0]} (T^*L \times V)/\Gamma $ to $V\subset \ker
d\Phi (x) \subset T_x M$ by the identity map.
Define $\th = \sigma ^* (h - \langle\Phi | \xi\rangle )$.  Because
$G$ is assumed to be abelian, $\xi$ is trivially in the center of
$\fg$, and so $\th$ is $G$-invariant.
At this juncture it is convenient to pretend that $U_0 = (T^*L \times
V)/\Gamma$.  By the invariance of $\th \in C^\infty ((T^*L \times
V)/\Gamma)^{L\times K}$, $\th ([l, \lambda, v]) = \th ([1, \lambda, v])$ for
all $(l,\lambda, v) \in L\times \fl^* \times V$.  Now define
$h_V (v) = \th ([1, 0, v])$.  This 
$h_V \in C^\infty (V)^{\Gamma\times K}$  is the desired
Hamiltonian: we have
$$
dh_V (0) = d\th([1, 0, 0])|_V = 0
$$ 
and 
$$d^2 h_V (0) = d^2\th ([1, 0, 0])|_V 
= d^2\left(\sigma^*(\augh)\right)([1, 0, 0])|_V 
= \sigma^*\left( d^2(\augh)(x)|_V\right).
$$
 
	Third and last, we prove that \rpo s for $h_V$ give rise to
\rpo s for $h$.  
Consider the fully reduced system $(P, h_P)$ obtained by reducing the
system $(V, h_V)$ by the action of $\Gamma\times K$.  We can arrive
at $(P, h_P)$ also by reducing the system 
$( (T^*L\times V)/\Gamma, \th\, )$ by the action of 
$G = L\times K$.  Thus, \rpo s for $h_V$ correspond to periodic
orbits for $h_P$, which in turn give rise to \rpo s for $\th$, or
equivalently to \rpo s for $\augh$.  But \rpo s for $\augh$ are
\rpo s for $h$.  Indeed, writing $\varphi_f^t$ for the Hamiltonian
flow of $f$, we have
$\varphi_{\augh}^t = 
\varphi_{-\langle\Phi | \xi \rangle}^t \circ \varphi_h^t$
by the $G$-invariance of $h$.  Therefore, 
$\varphi_{\augh}^t(x) = g\cdot x$ for some $g\in G$ if and only
if $\varphi_h^t(x) = g^{\prime}\cdot x$ with 
$g^{\prime} = \exp(-t\xi)\, g\in G$.
This concludes the proof of 
Theorem~\ref{Theorem2}.

\section{Concluding remarks}
\noindent
I.   Corollary proves the existence of \rpo s
for a pair of axially symmetric rigid bodies subject to a
coupling potential which depends on the angle between the bodies.  We
expect the theorem to prove the existence of \rpo s
for the motion of Riemann ellipsoids.  \\

\noindent
II.  In Theorem~\ref{Theorem2}, we compute $d^2h_V(0)$ from 
$d^2(\augh)(x)|_{\ker d\Phi(x)}$ by `killing' the direction of
degeneracy.  This computation needs only the {\it existence\/} of
the `Darboux coordinates' $\sigma : (T^*L \times V)/\Gamma \supset U_0 \to
M$ (see the proof above), but by computing without the explicit
form of $\sigma$, we lose some information.
Alternatively we can compute the Taylor expansion of $\sigma$
and translate the higher-order information on the jet of the original 
Hamiltonian $h$ into the information on the jet of the model 
Hamiltonian $h_V$.


\begin{thebibliography}{XXXX}


\bibitem[AMM]{AMM} J. Arms, J. Marsden, and V. Moncrief, Symmetry and
bifurcations of momentum mappings, {\em Comm.\ Math.\ Phys.} {\bf 78}
(1980/81), no. 4, 455--478.

\bibitem[B]{Ba} T. Bartsch, A generalization of the Weinstein-Moser
theorems on periodic orbits of a Hamiltonian system near an
equilibrium, {\em Ann.\ Inst.\ H. Poincar\'e Anal.\ Non Lin\'eaire\/}
{\bf 14} (1997), no. 6, 691--718.

\bibitem[GLS]{GLS} V. Guillemin, E. Lerman, and S. Sternberg, {\em
Symplectic fibrations and multiplicity diagrams}, Cambridge
University Press, Cambridge, 1996.

\bibitem[GS]{GS} V. Guillemin and S. Sternberg, A normal form for the
moment map, in: {\em Differential Geometric Methods in Mathematical
Physics\/} (ed. S. Sternberg), Reidel, Dordrecht, 1984.

\bibitem[L]{L} E. Lerman, On the centralizer of invariant functions on
a Hamiltonian $G$-space, {\em J. Diff.\ Geom.} {\bf 30} (1989)
805--815.

\bibitem[LS]{LS} E. Lerman and S. F. Singer, Stability and persistence
of relative equilibria at singular values of the moment map,
{\em Nonlinearity\/} {\bf 11} (1998) 1--13.

\bibitem[Mr]{Mr} C.-M. Marle, Mod\`ele d'action hamiltonienne d'un groupe
de Lie sur une vari\'et\'e symplectique, {\em Rend. del Seminario
Matematico\/} {\bf 43} (1985) 227--251.

\bibitem[MnRS]{MnRS} J. A. Montaldi, R. M. Roberts, and I. N. Stewart,
Periodic solutions near equilibria of symmetric Hamiltonian systems,
{\em Phil. Trans. R. Soc. Lond.} A {\bf 325} (1988) 237--293.

\bibitem[Ms]{Ms} J. Moser, Periodic orbits near an equilibrium 
   and a theorem of A. Weinstein, {\em Pure Appl. Math.} {\bf 29}
   (1976) 727--747.

\bibitem[O]{Ortega} J.-P. Ortega, {\em Symmetry, reduction and
stability in Hamiltonian systems}, thesis, University of California,
Santa Cruz, 1998.

\bibitem[OR]{OR} J.-P. Ortega and T. Ratiu, Stability of Hamiltonian
relative equilibria, preprint, 1997.

\bibitem[SL]{SL} R. Sjamaar and E. Lerman, Stratified symplectic spaces
and reduction, {\em Ann.\ of Math.} {\bf 134 }(1991) 375--422.

\bibitem[W1]{W1} A. Weinstein, Normal modes for nonlinear Hamiltonian
   systems, {\em Invent. Math.} {\bf 20} (1973) 47--57.

\bibitem[W2]{W2} A. Weinstein, Bifurcations and Hamilton's principle, 
   {\em Math. Z.} {\bf 159} (1978) 235--248.

\end{thebibliography}
\end{document}